\documentclass[11pt]{article} 
\usepackage[fleqn]{amsmath}                  
\usepackage{amssymb}                  
\usepackage{amsthm}                   
\usepackage[]{geometry}                 
\usepackage{graphicx}                 
\usepackage{float}                    
\usepackage{caption}                  
\usepackage{cite}
\usepackage{blkarray}
\usepackage{booktabs}   
\usepackage{array}      
\usepackage{multirow}   
\usepackage{breqn}

\theoremstyle{definition}

\newcommand{\nn}{\nonumber}

\begin{document}
\title{A discrete-time $GI^{X}/Geo^{Y}/1$ queue with early arrival system}
\author{U. C. Gupta, F. P. Barbhuiya, Arunava Maity\thanks{E-mail addresses: umesh@maths.iitkgp.ac.in(U. C.
Gupta), faridaparvezb@gmail.com(F. P.
Barbhuiya), arunavamath@gmail.com(Arunava Maity).}\\
{\it Department of Mathematics, Indian Institute of Technology Kharagpur,}\\
{\it Kharagpur-721302, India}}
\maketitle{} \textbf{Abstract:} A discrete-time batch service queue
with batch renewal input and random serving capacity rule under the
late arrival delayed access system (LAS-DA), has recently appeared
in the literature \cite{barbhuiya2019discrete}. In this paper, we
consider the same model under the early arrival system (EAS), since
it is more applicable in telecommunication systems where an arriving
batch of packets needs to be transmitted in the same slot in which
it has arrived. In doing so, we derive the steady-state queue length
distributions at various epochs, and show that in limiting case the
result gets converted to the continuous-time queue
\cite{barbhuiya2019difference}. We discuss few numerical results as
well. \\
%
\section{Introduction}
\label{introd6}
%
%
The continuous-time batch-arrival and batch-service $GI^{X}/M^{Y}/1$
queue has been studied in the past by Economou and Fakinos
\cite{economou2003} and Cordeau and Chaudhry \cite{leo2009}.
Whereas, the former one gave the probability generating function
(pgf) of the system content distribution at points of arrivals, the
later one inverted the pgf given in [9] using roots method.
Recently, Barbhuiya and Gupta \cite{barbhuiya2019difference}
revisited the work in \cite{economou2003} and \cite{leo2009} and
proposed a different methodology for the analysis based on
difference equation technique, which is analytically as well as
computationally tractable.
\par Observing the limitations of continuous-time queue in
modeling digital communications, computer networks, tele-traffic
processes and many more, very recently, Barbhuiya and Gupta
\cite{barbhuiya2019discrete} considered the theoretical and
computational aspects of the discrete-time $GI^X/Geo^Y/1$ queue.
They studied the model under the assumption of only late arrival
with delayed access system (LAS-DA), according to which, an arrival
cannot depart in the same slot in which it has arrived even if the
server is idle, as opposed to the early arrival system (EAS).
Discrete-time queues with EAS policy are more significant in the
modeling of systems where the packets (information) needs to be
transmitted in the same slot in which it has arrived, under an
urgent situation. With an aim to fill the gap in the literature,
this paper studies the $GI^X/Geo^Y/1$ queue with EAS assumption. The
procedure used throughout the analysis runs parallel to that of
\cite{barbhuiya2019discrete}, however, the result differs
significantly.
\par The main contributions of the paper are as follows:
{(\em{i})} we obtain the steady-state queue-length distributions at
pre-arrival and arbitrary epochs in a readily presentable form.
Consequently, the numerical computation of state probabilities
becomes much easier and quite straight forward, {(\em{ii})} the
analysis provides an alternative methodology to solve many special
queueing models considered in the past, such as $GI^X/Geo/1$ (see
\cite{vinck1994analyzing}, \cite{chaudhry1997queue}), $GI/Geo/1$
(see \cite{chaudhry1996relations}), {(\em{iii})} we derive the
results of the continuous-time $GI^X/M^Y/1$ queue from the results
of the discrete-time counterpart.
\par The remaining portion of the paper is organized as follows. In
section \ref{md6} we provide the model description which is followed
by the mathematical analysis and the derivation of queue-length
distributions in section \ref{ge6}. In section \ref{sc6} we discuss
few special cases of the model. In section \ref{sec5} we obtain the
results of the continuous-time queue from the discrete-time
counterpart and finally present some numerical results in section
\ref{nr6} which is followed by the conclusion.
\section{Model description} \label{md6}
Though the model has been described in detail in
\cite{barbhuiya2019discrete}, for the sake of completeness, in this
section we give a brief overview of it.
\par Customers arrive in batches of random size $X$ with probability
mass function (pmf) $P(X=i)=g_i$, $i=1,2,3...,b$, probability
generating function (pgf) $G(z)=\sum_{i=1}^{b}g_{i}z^{i}, |z| \leq
1$ and the mean batch size $\overline{g}=\sum_{i=1}^{b}ig_{i}$. Here
$b \in \mathbb{N}$ is the maximum possible size of the arriving
batch. The inter-arrival times of the batches are independent and
identically distributed (i.i.d) random variables with common pmf
$a_n\!=\!P(A\!=\!n)$, $n\!\ge\! 1$~(where $A$ is the random variable
corresponding to the inter-arrival time), pgf
~$A(z)=\sum_{n=1}^{\infty}a_nz^n$, and mean inter-arrival time~$a$ =
$\frac{1}{\lambda}=A'(1)$. The customers are served in batches by a
single server, where the batch-size to be served is a random
variable ($Y$) with pmf $P(Y=i)=y_i$, $i=1,2,3...$, pgf
$Y(z)=\sum_{i=1}^{\infty}y_{i}z^{i}, |z| \leq 1$ and the mean
batch-size $\overline{y}$. At each service initiation stage the
server will serve minimum of \{$i$, whole queue\} with probability
$y_i$. The service times $S$ of the batches are independent and
geometrically distributed as $P(S=n)=b_n=(1-\mu)^{n-1}\mu,~ 0 < \mu
< 1,~ n\ge 1$, with mean $\frac{1}{\mu}$. The arrival and service
processes are independent of each other.  The traffic intensity
$\rho = \frac{\lambda \overline{g}}{\mu \overline{y}}$ and $\rho <
1$ ensures the stability of the system.

\section{Analysis of the model}\label{ge6}
\begin{figure}[h!]
\centering
  \includegraphics[width=15cm, height=7cm]{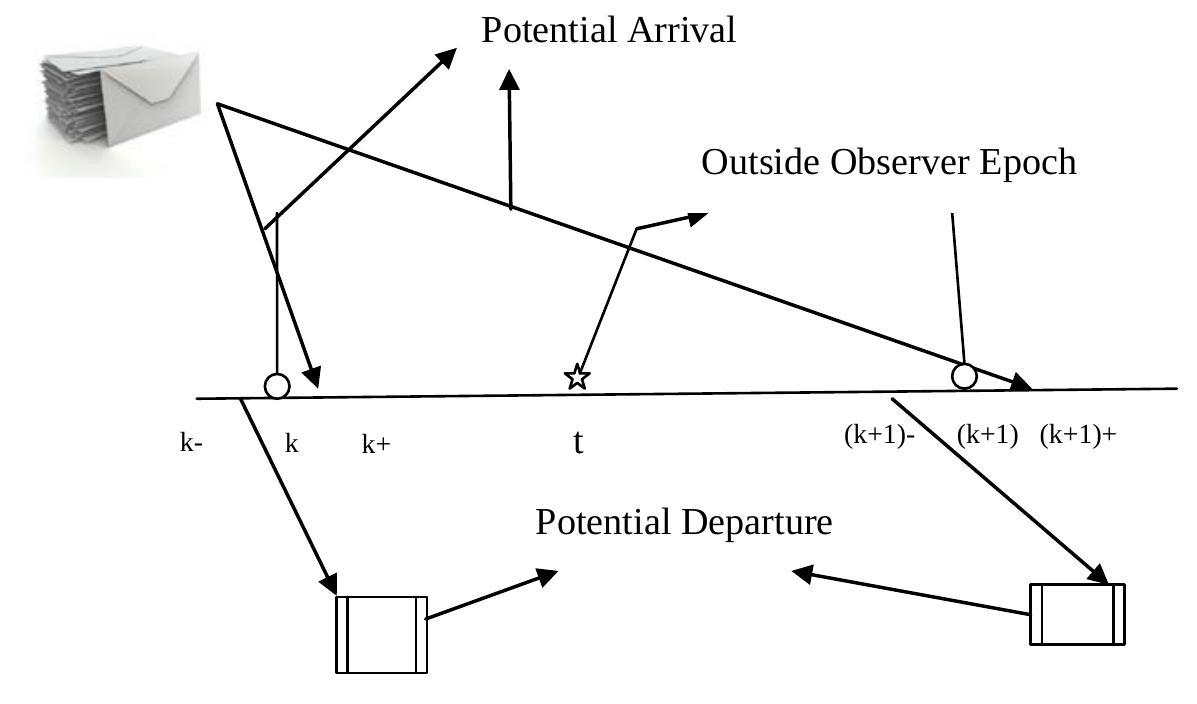}
  \caption{Various epochs for Early Arrival System (EAS)}\label{fig6.1}
\end{figure}
We consider the early arrival system (EAS) for the present model and
accordingly, let the time axis be divided into slots of equal length
such that the slot boundaries are marked by $ 0,1,2,....,k,...$. We
assume that the potential batch arrival occurs in $(k,k+$) and the
potential batch departure occurs in $(k-,~k)$ (see
Fig.\ref{fig6.1}).
The state of the system just before a potential batch arrival, i.e.,
at the instant $k$ is described by two random variables $(N_k, U_k)$
where $N_k$ is the number of customers in the queue and $U_k$ is the
remaining inter-arrival time of the next batch at the instant $k$.
We define the joint distribution of  $N_k$ and $U_k$  as,
$$\widehat{p}_n(k,u)= P[N_k=n, U_k=u], ~ u \geq 0, n \geq 0. $$
Relating the states of the system at two consecutive time epochs $k$
and $k+1$, using the arguments of supplementary variable technique
(SVT) by assuming the remaining inter-arrival time of the next batch
as the supplementary variable ($u \geq 1$), we obtain the set of
governing equations as,
\begin{eqnarray}
\widehat{p}_0(u-1, k+1) &=& \widehat{p}_0(u,k) + \mu
\sum_{i=1}^{\infty}\widehat{p}_i(u,k)\sum_{j=i}^{\infty}y_j + a_u
\mu
\sum_{i=0}^{\infty}\widehat{p}_i(0,k)\sum_{m=1}^bg_m\sum_{j=i+m}^{\infty}y_j,
\label{91} \\
\widehat{p}_n(u-1,k+1)&=& \widehat{p}_n(u,k)(1- \mu)+ \mu
\sum_{i=1}^{\infty}\widehat{p}_{n+i}(u,k)y_i + a_u(1-
\mu)\sum_{i=1}^ng_i\widehat{p}_{n-i}(0,k) \nn \\ && + a_u \mu
\sum_{i=1}^{\infty}y_i\sum_{m=1}^{min\{i+n, b\}}g_m
\widehat{p}_{n-m+i}(0,k),~~1 \leq n \leq b-1, \label{92}\\
\widehat{p}_n(u-1, k+1)&=& \widehat{p}_n(u,k)(1- \mu) + \mu
\sum_{i=1}^{\infty}\widehat{p}_{n+i}(u,k)y_i +
a_u(1-\mu)\sum_{i=1}^bg_i\widehat{p}_{n-i}(0,k) \nn \\ && + a_u \mu
\sum_{i=1}^{\infty}y_i\sum_{m=1}^{b}g_m \widehat{p}_{n-m+i}(0,k),~~n
\geq b. \label{93}
\end{eqnarray}
Further, in steady-state we define $ p_n(u)=\lim_{k \rightarrow
\infty}\widehat{p}_n(k,u) $ and thus (\ref{91})-(\ref{93}) reduces
to
\begin{eqnarray}
p_0(u-1) &=& p_0(u) + \mu
\sum_{i=1}^{\infty}p_i(u)\sum_{j=i}^{\infty}y_j + a_u \mu
\sum_{i=0}^{\infty}p_i(0)\sum_{m=1}^bg_m\sum_{j=i+m}^{\infty}y_j.
\label{1} \\
p_n(u-1)&=& p_n(u)(1- \mu)+ \mu \sum_{i=1}^{\infty}p_{n+i}(u)y_i +
a_u(1- \mu)\sum_{i=1}^ng_ip_{n-i}(0) \nonumber \\ && + a_u \mu
\sum_{i=1}^{\infty}y_i\sum_{m=1}^{min\{i+n, b\}}g_m p_{n-m+i}(0), ~~~ 1
\leq n \leq b-1. \label{2}\\
p_n(u-1)&=& p_n(u)(1- \mu)+ \mu \sum_{i=1}^{\infty}p_{n+i}(u)y_i +
a_u(1- \mu)\sum_{i=1}^bg_ip_{n-i}(0) \nonumber \\ && + a_u \mu
\sum_{i=1}^{\infty}y_i\sum_{m=1}^{b}g_m p_{n-m+i}(0), ~~~ n \geq b.
\label{3}
\end{eqnarray}
In order to obtain the steady-state queue-length distributions we
consider the transform $p_n^*(z)= \sum_{u=0}^{\infty}p_n(u)z^u$ and
$A(z)= \sum_{u=0}^{\infty}a_uz^u$ such that $a_0=0$ and $|z| \leq
1$. This gives $ p_n^*(1)=\sum_{u=0}^{\infty}p_n(u)\equiv p_n$.
Thus, using the transform as defined above, we obtain from
(\ref{1})-(\ref{3}) the following
\begin{eqnarray}
(z-1)p_0^*(z)&=& \mu
\sum_{i=1}^{\infty}p_i^*(z)\sum_{j=i}^{\infty}y_j -p_0(0) -\mu
\sum_{i=1}^{\infty}p_i(0)\sum_{j=i}^{\infty}y_j \nonumber \\ && +\mu
A(z)\sum_{i=0}^{\infty}p_i(0)\sum_{m=1}^bg_m\sum_{j=i+m}^{\infty}y_j, \label{5}\\
\left[z-(1- \mu)\right]p_n^*(z) &=& \mu
\sum_{i=1}^{\infty}y_ip_{n+i}^*(z)- (1- \mu)p_n(0) -
\mu\sum_{i=1}^{\infty}y_ip_{n+i}(0)\nonumber \\ && + A(z)(1- \mu)
\sum_{i=1}^ng_ip_{n-i}(0) \nn \\ && + \mu
A(z)\sum_{i=1}^{\infty}y_i\sum_{m=1}^{min\{i+n,
b\}}g_mp_{n-m+i}(0),~~1 \leq n \leq b-1. \label{6}\\
\left[z-(1- \mu)\right]p_n^*(z) &=& \mu
\sum_{i=1}^{\infty}y_ip_{n+i}^*(z)- (1- \mu)p_n(0) -
\mu\sum_{i=1}^{\infty}y_ip_{n+i}(0) \nn \\ && + A(z)(1- \mu)
\sum_{i=1}^bg_ip_{n-i}(0) + \mu
A(z)\sum_{i=1}^{\infty}y_i\sum_{m=1}^{b}g_mp_{n-m+i}(0),n \geq b.
\label{7}
\end{eqnarray}
Adding (\ref{5})-(\ref{7}), taking limit as $z\rightarrow 1$ and
using the normalizing condition $\sum_{n=0}^{\infty}p_n=1$ we obtain
\begin{eqnarray}
\sum_{n=0}^{\infty} p_n(0)=\frac{1}{a_1}=\lambda. \label{8}
\end{eqnarray}
The above result may also be interpreted intuitively as, $p_n(0)$
represents the mean number of times the remaining inter-arrival time
hits $zero$ per unit time, the sum of which for all $n$ becomes the
arrival rate $\lambda$. Let us now define $p_{n}^{-}$ as the
probability that the queue-length is $n$ just before the arrival of
a batch, i.e., at pre-arrival epoch. As $p_{n}^{-}$ is proportional
to $p_n(0)$ and $\sum_{n=0}^{\infty}p_{n}^{-}=1$, we have the
relation
\begin{eqnarray}\label{9}
p_n^-=\frac{p_n(0)}{\sum_{k=0}^{\infty}p_k(0)}=\frac{p_n(0)}{\lambda}.
\end{eqnarray}
%
%
\subsection{Steady-state distribution at pre-arrival and arbitrary
epochs}\label{dist6}
For the sequence of probabilities $\{p_n(0)\}_{0}^{\infty}$ and
$\{p_n^*(z)\}_{0}^{\infty}$ we define the right shift operator $D$
as $Dp_n(0)=p_{n+1}(0)$ and $Dp_n^*(z)=p_{n+1}^*(z)$ for all $n$.
Thus (\ref{7}) can be rewritten as
\begin{eqnarray}
\left[z-(1- \mu + \mu \sum_{i=1}^{\infty}y_iD^i)\right]p_n^*(z)&=&
\left[-(1-\mu)D^b-\mu\sum_{i=1}^{\infty}y_iD^{i+b}+ A(z)(1-\mu)
\sum_{i=1}^bg_iD^{b-i}\right. \nn \\ &&\left.+A(z)\mu
\sum_{m=1}^bg_m\sum_{i=1}^{\infty}y_iD^{b+i-m}\right]p_{n-b}(0),~~
n\geq b.\label{10}
\end{eqnarray}
Substituting $z=1- \mu + \mu \sum_{i=1}^{\infty}y_iD^i$ in
(\ref{10}) we get
\begin{eqnarray}
\left[-(1-\mu)D^b-\mu\sum_{i=1}^{\infty}y_iD^{i+b}+ A(1- \mu + \mu
\sum_{i=1}^{\infty}y_iD^i)(1-\mu) \sum_{i=1}^bg_iD^{b-i}\right. && ~
\nn \\\left.+ A(1- \mu + \mu
\sum_{i=1}^{\infty}y_iD^i)\mu\sum_{m=1}^bg_m\sum_{i=1}^{\infty}y_iD^{b+i-m}
\right]p_{n}(0)&=&0,~~{n \geq 0},\label{11}
\end{eqnarray}
which is a homogeneous difference equation with constant coefficient
with the corresponding characteristic equation (c.e.) as $\left(
A(1- \mu + \mu Y(s))\sum_{i=1}^bg_is^{b-i}-s^b\right)\left(1- \mu +
\mu Y(s)\right)=0$. But since
\begin{eqnarray}
A(1- \mu + \mu Y(s))\sum_{i=1}^bg_is^{b-i}-s^b=0 \label{12}
\end{eqnarray}
has exactly $b$ roots inside the unit circle under the stability
condition $\rho = \frac{\lambda \overline{g}}{\mu \overline{y}} < 1$
(see Barbhuiya and Gupta \cite{barbhuiya2019discrete}), hence the
c.e. has exactly $b$ roots inside $|s|=1$ denoted by $r_1, r_2, . .
., r_b$. Thus the general solution of (\ref{11}) is given by
\begin{eqnarray}\label{13}
p_{n}(0)&=& \sum_{i=1}^{b}c_i r_i^n~,~n \ge 0,
\end{eqnarray}
where $c_1, c_2, . . ., c_b$ are the arbitrary constants yet to be
determined. Substituting the expression of $p_n(0)$ from (\ref{13})
in (\ref{10}) we obtain
\begin{eqnarray}
\left(z-(1- \mu + \mu \sum_{i=1}^{\infty}y_iD^i)\right)p_n^*(z)&=&
-(1-\mu)\sum_{j=1}^bc_jr_j^{n}-\mu\sum_{i=1}^{\infty}y_i\sum_{j=1}^bc_jr_j^{n+i}
\nn \\ && + A(z)(1-\mu) \sum_{i=1}^bg_i\sum_{j=1}^bc_jr_j^{n-i} \nn
\\ && +A(z)\mu
\sum_{m=1}^bg_m\sum_{i=1}^{\infty}y_i\sum_{j=1}^bc_jr_j^{n+i-m},~~
{n \geq b}.\label{614}
\end{eqnarray}
which is a non-homogeneous difference equation and its general
solution is given by
\begin{eqnarray}
p_n^*(z)= \sum\limits_{j}^{}d_js_j^n(z) + \sum_{j=1}^bc_j\left\{
\frac{(1- \mu + \mu
Y(r_j))(A(z)\sum\limits_{i=1}^bg_ir_j^{b-i}-r_j^b)}{z-1 +\mu - \mu
Y(r_j)} \right\}r_j^{n-b},~{n \geq b}.\label{18}
\end{eqnarray}
In (\ref{18}) the first term in the R.H.S is the solution
corresponding to the homogeneous part of (\ref{614}), the second
term is the particular solution, $s_j(z)$'s are the roots of $z-1
+\mu - \mu Y(s)=0$ for a fixed $z$ and $d_j$'s are the corresponding
arbitrary constants. Since $\sum_{n=0}^{\infty}p_{n}=1$, so
$\sum_{n=b}^{\infty}p_{n} \leq 1$ i.e.,
$\sum_{n=b}^{\infty}p_{n}^{*}(1) \leq 1$. Setting $z=1$ and summing
over $n$ from $b$ to $\infty$ in (\ref{18}), we have the first term
in R.H.S as $\sum\limits_jd_j(\sum\limits_{n=b}^{\infty}s_j^n(1))$,
where $s_j(1)$ are the roots of $Y(s)=1$. These roots lie outside
and on the unit circle $|s|=1$ (for proof, see Barbhuiya and Gupta
\cite{barbhuiya2019discrete}) and hence
$\sum_{n=b}^{\infty}s_j^n(1)$ diverges. Thus to ensure the
convergence of (\ref{18}) we must have $d_j=0$ for all $j$ and
consequently we get the general solution of (\ref{614}) as
\begin{eqnarray}
p_n^*(z)= \sum_{j=1}^bc_j\left\{ \frac{(1- \mu + \mu
Y(r_j))(A(z)\sum\limits_{i=1}^bg_ir_j^{b-i}-r_j^b)}{z-1 +\mu - \mu
Y(r_j)} \right\}r_j^{n-b},~{n \geq b}.\label{19}
\end{eqnarray}
For $p_n^*(z)$, $1 \leq n \leq b-1$, we seek a similar expression as
given in (\ref{19}). Thus we need to find the conditions under which
$p_n^*(z)$, $1 \leq n \leq b-1$ satisfies (\ref{6}). Substituting
the respective expressions in (\ref{6}) we obtain
\begin{eqnarray}
\sum_{j=1}^bc_jr_j^n\left\{ (1- \mu)\sum_{i=n+1}^bg_ir_j^{-i}+ \mu
\sum_{i=1}^{\infty}y_ir_j^i\left(\sum_{m=1}^bg_mr_j^{-m} -
\sum_{m=1}^{min(i+n,b)}g_mr_j^{-m}\right) \right\} &=& 0, \nn \\ {1
\leq n \leq b-1}. &&  \label {20}
\end{eqnarray}
Putting $n=b-1,~b-2,...,1$ in (\ref{20}) and using the condition
that $g_b\neq 0$, we obtain the following set of $b-1$ equations:
\begin{eqnarray} \label{21}
\sum_{j=1}^{b} \frac{c_j}{r_j}=\sum_{j=1}^{b}
\frac{c_j}{r_j^2}=...=\sum_{j=1}^{b}
\frac{c_j}{r_j^{b-2}}=\sum_{j=1}^{b} \frac{c_j}{r_j^{b-1}}=0.
\end{eqnarray}
%
Also summing over $n$ from $0$ to $\infty$ in (\ref{13}) and using
(\ref{8}) we obtain
\begin{eqnarray} \label{22}
\lambda = \sum_{i=1}^{b}\frac{c_i}{1-r_{i}}.
\end{eqnarray}
Solving the system of $b$ equations (\ref{21}) and (\ref{22}), we
can obtain the constants $c_j$ for $j=1,2,...,b$. This makes the
expression of $p_n(0)$ ($n \geq 0$) given in (\ref{13}) completely
known. Moreover, $p_n^*(z)$ is given by
\begin{eqnarray}
p_n^*(z)= \sum_{j=1}^bc_j\left\{ \frac{(1- \mu + \mu
Y(r_j))(A(z)\sum\limits_{i=1}^bg_ir_j^{b-i}-r_j^b)}{z-1 +\mu - \mu
Y(r_j)} \right\}r_j^{n-b},~n \geq 1.\label{23}
\end{eqnarray}
Therefore, using (\ref{9}), (\ref{13}) and (\ref{23}) we get the
explicit expression of the system content distributions at
pre-arrival ($p_n^-$) and arbitrary ($p_n$) epochs as
\begin{eqnarray}
p_n^- &=& \frac{1}{\lambda}\sum\limits_{i=1}^{b} c_i r_i^n ~,~n \geq
0, \label{24} \\
p_n &=& p_n^*(1) = \sum_{j=1}^bc_j\left\{ \frac{(1- \mu + \mu
Y(r_j))(\sum\limits_{i=1}^bg_ir_j^{b-i}-r_j^b)}{\mu (1- Y(r_j))}
\right\}r_j^{n-b},~n \geq 1, \label{25}\\
p_0 &=& 1- \sum_{n=1}^{\infty}p_n = 1- \sum_{j=1}^bc_j\left\{
\frac{(1- \mu + \mu
Y(r_j))(\sum\limits_{i=1}^bg_ir_j^{1-i}-r_j)}{\mu (1-r_j) (1-
Y(r_j))} \right\}. \label{26}
\end{eqnarray}
This completes the analysis of $GI^X/Geo^Y/1$ queue under EAS
policy. The results of many special cases of this model can be
obtained directly from (\ref{24})-(\ref{26}) as discussed in the
forthcoming section. However, we now provide a small example to
give an overview of the present analysis to the readers.\\\\
\textbf{Example:} Consider a discrete-time $Geo^X/Geo^Y/1$ queue,
i.e., the inter-arrival time follows geometric distribution. Suppose
$\lambda=0.2$, $\mu= 0.5$ and the mean inter-arrival time is
$a=\frac{1}{\lambda}$. Let the arriving batch size distribution is
$g_1=0.4$, $g_2=0.3$ and $g_3=0.3$ and service batch size
distribution is $y_1=0.4$ and $y_2=0.6$. This gives $b=3$, $G(z)=
0.4z+0.3z^2+0.3z^3$, $Y(z)=0.4z+0.6z^2$, $\overline{g}= 1.9$,
$\overline{y}=1.6$, $\rho= 0.475$ and $A(z)= \frac{\lambda
z}{1-(1-\lambda)z}$.
\par The root equation as deduced from equation (\ref{12}) is
\begin{eqnarray}
1- A(1-\mu+ \mu Y(s))\sum_{i=1}^3g_is^{-i} &=&0 \nonumber
\\ \textrm{i.e.,} ~~~0.24s^5+0.184s^4-0.566s^3+0.07s^2+0.042s+0.03 &=& 0.
\end{eqnarray}
The corresponding roots are $-2.002850$,
$-0.183817\pm0.263763i$,
$0.603819$, $1.0$. Clearly, $r_1=-0.183817+0.263763i$, $r_2=
-0.183817-0.263763i$ and $r_3=0.603819$ as $|r_i|<1$ for $i=1,2,3$.
Now the system of three equations deduced from (\ref{21}) and
(\ref{22}) are
\begin{eqnarray}
\frac{c_1}{-0.183817+0.263763i} + \frac{c_2}{-0.183817-0.263763i}+
\frac{c_3}{0.603819}&=&0 \nonumber \\
\frac{c_1}{(-0.183817+0.263763i)^2} +
\frac{c_2}{(-0.183817-0.263763i)^2}+
\frac{c_3}{(0.603819)^2}&=&0 \nonumber \\
\frac{c_1}{1.183817-0.263763i} + \frac{c_2}{1.183817+0.263763i} +
\frac{c_3}{0.396181} &=& 0.2 \nonumber
\end{eqnarray}
which can be solved to obtain $c_1= 0.027481-0.000834i$, $c_2=
0.027481+0.000834i$ and $c_3= 0.061593$. With the known values of
$c_i$'s and $r_i$'s for $i=1,2,3$, the steady state probabilities at
pre-arrival ($p_n^-$) and arbitrary ($p_n$) epochs can be directly
evaluated using (\ref{24}), (\ref{25}) and (\ref{26}) as $p_0^-=p_0=
0.582779$, $p_1^-=p_1=0.137643$, $p_2^-=p_2=0.101641$,
$p_3^-=p_3=0.076705$ and so on. Here both the state probabilities
are equal due to Bernoulli arrivals. Furthermore, the average
system-content at pre-arrival ($L^-$) and arbitrary ($L$) epochs are
given by $L^-=\sum_{n=1}^{N}np_n^-= 1.133649$ and
$L=\sum_{n=1}^{N}np_n= 1.133649$ for a sufficiently large $N(=500)$.
Meanwhile, as the value of $n$ becomes larger (say $n>30$), the
ratio $p_{n+1}^-/p_n^-$ tends to $0.603819$ which is $r_3$, the
unique largest real root inside the unit circle (see Barbhuiya and
Gupta \cite{barbhuiya2019discrete}).
\section{Special Cases}\label{sc6}
In this section, we discuss some special queueing models whose
results can be deduced directly from the analysis done in Section
\ref{dist6}.
\begin{itemize}
\item The model $GI^X/Geo^Y/1$ reduces to $GI^X/Geo/1$ if we consider
$y_1=1$ and $y_i=0$ for $i\geq 2$, i.e., the pgf $Y(s)=s$.
Accordingly, the system content distribution at pre-arrival and
arbitrary epochs reduces to
\begin{eqnarray}
p_n^- &=& \frac{1}{\lambda}\sum\limits_{i=1}^{b} c_i r_i^n ~,~n \geq
0, \nonumber \\
p_n &=& \sum_{j=1}^bc_j\left\{ \frac{(1- \mu + \mu
r_j)(\sum\limits_{i=1}^bg_ir_j^{b-i}-r_j^b)}{\mu (1- r_j)}
\right\}r_j^{n-b},~n \geq 1, \nonumber \\
p_0 &=& 1- \sum_{j=1}^bc_j\left\{ \frac{(1- \mu + \mu
r_j)(\sum\limits_{i=1}^bg_ir_j^{1-i}-r_j)}{\mu (1-r_j)^2} \right\}.
\nonumber
\end{eqnarray}
where $r_j$'s are the $b$ roots of $ A(1- \mu + \mu
s)\sum_{i=1}^bg_is^{b-i}-s^b=0$ inside the unit circle and the
constants $c_j$, $j=1,2,...b$ can be obtained by solving the system
of equations (\ref{21}) and (\ref{22}). This model was earlier
considered by Chaudhry and Gupta \cite{chaudhry1997queue} and Vinck
and Bruneel \cite{vinck1994analyzing}. However, the analysis carried
out in this paper provides and alternative approach to the solution
of the model.
\item Similarly, we can obtain the results for $GI/Geo^Y/1$ queue by
taking $g_1=1$ and $g_i=0$ for $i\geq 2$, i.e., the pgf $G(s)=s$.
The single root inside the unit circle corresponding to (\ref{12})
is denoted by $r_1$ and the constant obtained from (\ref{22}) is
$c_1= \lambda (1-r_1)$. Thus, the distribution at pre-arrival and
arbitrary epochs reduces to
\begin{eqnarray}
p_n^-&=& (1-r_1)r_1^n~, ~n \geq 0, \nonumber \\
p_n &=& \frac{\lambda (1-r_1)^2(1- \mu + \mu Y(r_1)) r_1^{n-1}}{\mu
(1- Y(r_1))}, ~ n\geq
1, \nonumber \\
p_0 &=& 1- \frac{\lambda (1- r_1)(1- \mu + \mu Y(r_1))}{\mu (1-
Y(r_1))}. \nonumber
\end{eqnarray}
\item Along the same line if we consider $g_1=1$, $g_i=0$ for $i\geq 2$ and $y_1=1$,
$y_i=0$ for $i\geq 2$, then $GI^X/Geo^Y/1$ queue reduces to the
$GI/Geo/1$ queue. Here the pgf's $G(s)=s$ and $Y(s)=s$, the single
root inside the unit circle corresponding to (\ref{12}) is denoted
by $r_1$ with the constant obtained as $c_1= \lambda (1-r_1)$.
Consequently, the probability distribution obtained at pre-arrival
and arbitrary epochs are given by
\begin{eqnarray}
p_n^-&=& (1-r_1)r_1^n~, ~n \geq 0, \nonumber \\
p_n &=& \frac{\lambda}{\mu}(1-r_1)(1- \mu + \mu r_1)r_1^{n-1},~ n \geq 1,\nonumber \\
p_0 &=& 1- \frac{\lambda}{\mu} (1- \mu + \mu r_1). \nonumber
\end{eqnarray}
These results exactly matches with the one obtained by Chaudhry et
al. \cite{chaudhry1996relations}.
\end{itemize}
\section{Results of continuous-time $GI^X/M^Y/1$ queue}\label{sec5}
In this section we convert the analytical results obtained for the
discrete-time queue to the corresponding continuous-time
$GI^X/M^Y/1$ queue. However, the results of discrete-time queues
cannot be derived from the continuous time counterpart (see Takagi
\cite{takagi1993queueing}). So for the continuous time case, we
assume the random variable $\widehat{A}$ as the inter-arrival times
between two successive batches which are i.i.d with distribution
function $\widehat{A}(x)$, Laplace-Stieltjes transform (L.S.T)
$A^*(s)$ and mean $\frac{1}{\widehat{\lambda}}$, where
$\widehat{\lambda}$ is the arrival rate of the batches. The service
times $\widehat{S}$ of the batches are exponentially distributed
with rate $\widehat{\mu}$. Let the time axis be divided into
intervals of equal length $\Delta$, such that $\Delta>0$ and is
sufficiently small. Thus we have, $a_n= Pr((n-1)\Delta < \widehat{A}
\leq n\Delta)$, $n \geq 1$ and $E(A)\Delta= E(\widehat{A})$. One may
note that $\widehat{\lambda}\Delta=\lambda$ and
$\widehat{\mu}\Delta=\mu$. We now prove that the c.e. (\ref{12})
reduces to the c.e. corresponding to the continuous-time queue
\cite{barbhuiya2019difference}.
\begin{eqnarray}
A(1- \mu + \mu Y(s)) &\equiv& \lim_{\Delta \rightarrow 0}A(1-
\widehat{\mu}\Delta + \widehat{\mu}\Delta Y(s))\nn \\ ~ &=&
\lim_{\Delta \rightarrow 0} \sum_{n=1}^{\infty}Pr((n-1)\Delta <
\widehat{A} \leq n\Delta)[1- \widehat{\mu}\Delta +
\widehat{\mu}\Delta Y(s)]^n \nn \\ ~ &=& \lim_{\Delta \rightarrow 0}
\sum_{n=1}^{\infty}
[\widehat{A}(n\Delta)-\widehat{A}((n-1)\Delta)]\left[
1-\frac{\widehat{\mu}(1-Y(s))\Delta n}{n} \right]^n \nn \\ ~ &=&
\int_{0}^{\infty}e^{-\widehat{\mu}(1-Y(s))x}d\widehat{A}(x)\nn \\
&=& A^*(\widehat{\mu}(1-Y(s))) \nn
\end{eqnarray}
Thus, $A^*(\widehat{\mu}(1-Y(s)))\sum_{i=1}^bg_is^{b-i}-s^b=0$
represents the c.e. corresponding to the continuous time queue,
which also have exactly $b$ roots ($\widehat{r}_1$, $\widehat{r}_2$,
..., $\widehat{r}_b$) inside the unit circle $|s|=1$. Now replacing
$r_j$'s by $\widehat{r}_j$'s, $\lambda$ by $\widehat{\lambda}\Delta$
in (\ref{21}) and (\ref{22}) and solving the new system of
equations, we obtain the constants in terms of $\Delta$, or to be
precise, as linear functions of $\Delta$ i.e., $\widehat{c}_j\Delta$
such that $\widehat{c}_j$ for $j=1,2,...,b$ are known. Again using
$\widehat{\lambda}\Delta=\lambda$, $\widehat{\mu}\Delta=\mu$,
$r_j=\widehat{r}_j$, $c_j=\widehat{c}_j\Delta$ and taking limit as
$\Delta \rightarrow 0$ in (\ref{24})-(\ref{26}) we obtain the
pre-arrival and arbitrary epoch probabilities as given in Barbhuiya
and Gupta \cite{barbhuiya2019difference}.
\begin{eqnarray}
p_n^- &=&
\frac{1}{\widehat{\lambda}}\sum_{j=1}^b\widehat{c}_j\widehat{r}_j^n,
~~n \geq 0 \nn \\
p_n &=& \sum_{j=1}^b\widehat{c}_j\left\{
\frac{G(\widehat{r}_j^{-1})-1}{\widehat{\mu}(1-
Y(\widehat{r}_j))}\right\}\widehat{r}_j^n, ~~n \geq 1 \nn \\
p_0 &=& 1- \sum_{j=1}^b\widehat{c}_j\left\{
\frac{\sum\limits_{i=1}^bg_i\widehat{r}_j^{1-i}-\widehat{r}_j}{\widehat{\mu}
(1-\widehat{r}_j) (1- Y(\widehat{r}_j))} \right\}. \nn
\end{eqnarray}
\section{Numerical Results}\label{nr6}
In this section, we present few numerical results in the self explanatory tables in order to
demonstrate an overview of the theoretical work done so far.
The numerical results are provided
correct upto six decimal places. The input parameters taken for
deterministic ($D$) inter-arrival time distribution in Table
\ref{table1} are same as given in
 Chaudhry and Gupta \cite{chaudhry1997queue}. It can be
observed that the numerical  values  for $p_n^-$ and $p_n$ exactly
matches with that of Chaudhry and Gupta \cite{chaudhry1997queue}, as
expected. In addition to that, we have considered arbitrary,
geometric ($Geo$) and deterministic ($D$) inter-arrival time
distributions, as given in Table \ref{table1} and \ref{table2}
respectively. Whereas, in Table \ref{table2} one may observe that
the service batch-size distribution is taken as geometric with
parameter 0.4, whose pgf is $Y(z)= \frac{0.4z}{1- 0.6z}$. A
noticeable feature from the distributions taken here is that, the
server can comply with infinite support also, and the service
batch-size distribution (with infinite support) need not to be
truncated to obtain the set of probabilities $p_n^-$ and $p_n$. It
is clearly an advantage  of our methodology developed here. Moreover
an important observation is that, in the limiting case the ratio of
the pre-arrival epoch probabilities $\frac{p_{n+1}^-}{p_n^-}$
converges to the largest root in absolute value of the c.e.
(\ref{12}) lying inside the unit circle (for analytical proof, one
may refer to Barbhuiya and Gupta \cite{barbhuiya2019discrete}). We
can thus infer that, the tail probabilities at pre-arrival epoch can
be approximated by the root of the c.e. with largest modulus inside
the unit circle.
\begin{table}[!htp]
\centering \caption{System content distribution at pre-arrival
($p_n^-$) and arbitrary ($p_n$) epochs}  \footnotesize
\label{table1}
\begin{tabular}{p{0.5cm}p{2cm}p{2cm}p{2cm}p{2cm}p{2cm}p{2cm}}
\hline ~ &\multicolumn {3}{c}{$GI\equiv D$} &
\multicolumn{3}{c}{$GI\equiv \textrm{arbitrary}$}\\
\cmidrule(lr){2-4}   \cmidrule(lr){5-7}  &
\multicolumn{3}{l}{$g_1=0.2, ~ g_5=0.3, ~g_{10}=0.5, ~\bar{g}=6.7$}
& \multicolumn{3}{l}{$ g_2=0.2,~ g_3=0.3, ~g_{5}=0.3, ~g_8=0.2,$}\\
& \multicolumn{3}{l}{$ b=10, ~Y(z)=z, ~ \mu =0.9, ~ \lambda =0.1,$}
& \multicolumn{3}{l}{$ \bar{g}=4.4, ~b=8, ~ \mu=0.4, ~\lambda
=0.1667, $}\\
& \multicolumn{3}{l}{$ \rho =0.7444$} & \multicolumn{3}{l}{$Y(z)=
0.3z+ 0.2z^3+0.25z^4+0.25z^6, $}\\
& \multicolumn{3}{l}{} & \multicolumn{3}{l}{$\rho =0.539, ~ a_7=0.5,
~a_{10}=0.2, ~a_{15}=0.3 $}\\
\cmidrule(lr){2-4}   \cmidrule(lr){5-7} $n$ & $p_n^-$ & $p_n$ &
$p_{n+1}^-/p_n^-$ & $p_n^-$ & $p_n$ &
$p_{n+1}^-/p_n^-$\\
\hline

0   &   0.578601    &   0.313415    &   0.253224    &   0.878471    &   0.406746    &   0.037114    \\
1   &   0.146516    &   0.067368    &   0.742023    &   0.032604    &   0.093348    &   0.804337    \\
2   &   0.108718    &   0.075652    &   0.613434    &   0.026224    &   0.130608    &   0.697163    \\
3   &   0.066691    &   0.081029    &   0.595849    &   0.018283    &   0.096706    &   0.892676    \\
4   &   0.039738    &   0.084169    &   0.599919    &   0.016321    &   0.067738    &   0.649924    \\
5   &   0.023840    &   0.068693    &   0.601121    &   0.010607    &   0.090820    &   0.465432    \\
6   &   0.014330    &   0.063502    &   0.600828    &   0.004937    &   0.016528    &   1.148267    \\
7   &   0.008610    &   0.060433    &   0.600749    &   0.005669    &   0.028922    &   0.711361    \\
8   &   0.005173    &   0.058478    &   0.600771    &   0.004033    &   0.049630    &   0.254162    \\
$\vdots$ & $\vdots$ &$\vdots$ & $\vdots$ & $\vdots$ & $\vdots$ & $\vdots$\\
124 &   0.000000    &   0.000000    &   0.600774    &   0.000000    &   0.000000    &   0.589450    \\
125 &   0.000000    &   0.000000    &   0.600774    &   0.000000    &   0.000000    &   0.589450    \\
126 &   0.000000    &   0.000000    &   0.600774    &   0.000000    &   0.000000    &   0.589450    \\
127 &   0.000000    &   0.000000    &   0.600774    &   0.000000    &   0.000000    &   0.589450    \\
$\vdots$ & $\vdots$ &$\vdots$ & $\vdots$ & $\vdots$ & $\vdots$ & $\vdots$\\

\hline
sum & 1.000000 & 1.000000 & ~ &  0.999999 & 0.999999 & ~ \\
mean & $L^- = 1.11$ & $L= 3.73$ & ~ & $L^-= 0.39$ &
$L= 2.27$ & ~ \\
\hline
\end{tabular}
\end{table}

\begin{table}[!htp]
\centering \caption{System-content distribution at pre-arrival
($p_n^-$) and arbitrary ($p_n$) epochs}  \footnotesize
\label{table2}
\begin{tabular}{p{0.5cm}p{2cm}p{2cm}p{2cm}p{2cm}p{2cm}p{2cm}}
\hline ~ &\multicolumn {3}{c}{$GI\equiv Geo$} &
\multicolumn{3}{c}{$GI\equiv D$}\\
\cmidrule(lr){2-4}   \cmidrule(lr){5-7}  &
\multicolumn{3}{l}{$g_2=0.2,~ g_5=0.25, ~g_{10}=0.3, ~g_{12}=0.25,$}
& \multicolumn{3}{l}{$g_1=0.3,~ g_4=0.3,~ g_6=0.4,~\bar{g}=3.9,~ b=6,$}\\
& \multicolumn{3}{l}{$\bar{g}=7.65,~ b=12,~ \mu=0.5,~
\lambda=0.2,~\rho=0.651$} & \multicolumn{3}{l}{$\mu =0.6,~ \lambda
=0.2,~ Y(z)= \frac{0.4z}{1- 0.6z},~ \rho =0.52$} \\
& \multicolumn{3}{l}{$ Y(z)= 0.2z^2+ 0.3z^4+0.4z^6+0.1z^7,$} &
\multicolumn{3}{l}{}\\
\cmidrule(lr){2-4}   \cmidrule(lr){5-7} $n$ & $p_n^-$ & $p_n$ &
$p_{n+1}^-/p_n^-$ & $p_n^-$ & $p_n$ &
$p_{n+1}^-/p_n^-$\\
\hline
0   &   0.359677    &   0.359677    &   0.063929    &   0.730813    &   0.480000    &   0.095772    \\
1   &   0.022994    &   0.022994    &   1.608902    &   0.069991    &   0.097442    &   0.855181    \\
2   &   0.036995    &   0.036995    &   0.627629    &   0.059855    &   0.077541    &   0.805875    \\
3   &   0.023219    &   0.023219    &   1.533170    &   0.048236    &   0.082722    &   0.723580    \\
4   &   0.035599    &   0.035599    &   0.935126    &   0.034902    &   0.086760    &   0.645443    \\
5   &   0.033289    &   0.033289    &   1.209728    &   0.022528    &   0.060251    &   0.592548    \\
6   &   0.040271    &   0.040271    &   0.448552    &   0.013349    &   0.059059    &   0.551151    \\
7   &   0.018064    &   0.018064    &   1.868597    &   0.007357    &   0.018567    &   0.660793    \\
8   &   0.033753    &   0.033753    &   0.524978    &   0.004862    &   0.013352    &   0.636454    \\
$\vdots$ & $\vdots$ &$\vdots$ & $\vdots$ & $\vdots$ & $\vdots$ & $\vdots$\\
124 &   0.000004    &   0.000004    &   0.926546    &   0.000000    &   0.000000    &   0.620481    \\
125 &   0.000004    &   0.000004    &   0.926586    &   0.000000    &   0.000000    &   0.620481    \\
126 &   0.000004    &   0.000004    &   0.926549    &   0.000000    &   0.000000    &   0.620481    \\
127 &   0.000004    &   0.000004    &   0.926583    &   0.000000    &   0.000000    &   0.620481    \\
$\vdots$ &$\vdots$ &$\vdots$ & $\vdots$ & $\vdots$ & $\vdots$ & $\vdots$ \\
\hline
sum & 1.000000 & 1.000000 & ~ &  0.999999 & 1.000000 & ~ \\
mean & $L^-= 9.63$ & $L= 9.63$ & ~ & $L^- = 0.84$ &
$L= 1.99$ & ~ \\
\hline
\end{tabular}
\end{table}

\newpage
\section{Conclusion}
\label{concl6} In this paper, we have studied a discrete time
$GI^X/Geo^Y/1$ queue under early arrival system (EAS). We have used
supplementary variable technique and difference equation approach
and obtained the steady-state queue-length distribution at
pre-arrival and arbitrary epochs in a tractable way. The notable
feature of this novel approach is that, it does not require the
formation of any complex transition probability matrix, that usually
appears in conventional methodologies. We have derived the
analytical results of the continuous-time queue from the
discrete-time counterpart. Also, we have discussed some numerical
results along with some special cases of our model that can be
directly derived from the analysis. The results thus obtained may be
used in discrete-time systems arising in telecommunications and
computer networks where emergency cases needs to be addressed during
the transmission of packets.

\bibliographystyle{plain}  
\bibliography{mybibfile}

\end{document}